\input amstex
\input epsf
\documentstyle{amsppt}
\NoBlackBoxes

\def\ls{\leqslant}
\def\gs{\geqslant}

\def\l{\ell}
\def\ii{\bold {i}}
\def\RR{\Bbb R}
\def\FF{\Bbb F}
\def\C{\Omega}
\def\OO{\Cal O}
\def\veps{\varepsilon}
\def\hh{\frak h}
\def\gg{\frak g}
\def\Norm{\operatorname{Norm}}

\TagsOnRight

\topmatter
\title
The number of connected components in the double Bruhat cells
for nonsimply-laced groups
\endtitle
\author Michael Gekhtman$^*$, Michael Shapiro$^\ddag$,
 and Alek Vainshtein$^\dag$ \endauthor
\affil $^*$ Department of Mathematics, University of Notre Dame,
Notre Dame, IN 46556
{\tt Michael.Gekhtman.1\@nd.edu}\\
$^\ddag$ Matematiska Institutionen, KTH, Stockholm,
{\tt mshapiro\@math.kth.se}\\
$^\dag$ Departments of Mathematics and of Computer Science,
University of Haifa,  Israel 31905,
{\tt alek\@mathcs.haifa.ac.il}
\endaffil

\abstract
We compute the number of connected components in a generic real double
Bruhat cell for series $B_n$ and $C_n$ and an exceptional group $F_4$.
\endabstract

\thanks{{\smc 2000 Mathematics Subject Classification}: Primary 20F55;
Secondary 05E15, 14M15}
\endthanks


\rightheadtext{Double Bruhat cells for nonsimply-laced groups}
\leftheadtext{Gekhtman, Shapiro, and Vainshtein}
\endtopmatter


\document
\heading 1. Introduction and Main Result \endheading

 Let $G$ be a simply connected
semisimple algebraic group.
Let $B$ and $B_-$ be two $\RR$-split
opposite Borel subgroups, $N$ and $N_-$ their unipotent radicals,
$H = B \cap B_-$ an $\RR$-split maximal torus of $G$, and $W =
{\Norm}_G (H) / H$ the Weyl group of $G$.

The group $G$ has two {\it Bruhat decompositions}, with respect
to $B$ and $B_-\,$:
$$
G = \bigcup_{u \in W} B u B = \bigcup_{v \in W} B_- v B_-  \ .
$$
The {\it double Bruhat cells\/} $G^{u,v}$ are
defined by $G^{u,v} = B u B  \cap B_- v B_- \,$.
The maximal torus $H$ acts freely on $G^{u,v}$ by left (or right)
translations. The quotient of $G^{u,v}$ by this action is called
the {\it reduced double Bruhat cell\/} $L^{u,v} \subset G^{u,v}$
(see \cite{SSVZ, Z} for a more rigorous definition). Thus,
$G^{u,v}$ is biregularly isomorphic to $H \times L^{u,v}$, and all
properties of $G^{u,v}$ can be translated in a straightforward way
into the corresponding properties of $L^{u,v}$ (and vice versa).
In particular, Theorem~1.1 in \cite{FZ} implies that $L^{u,v}$ is
biregularly isomorphic to a Zariski open subset of an affine space.

Let $L^{u,v}(\RR)$ denote the real part of $L^{u,v}$, that is,
$L^{u,v}(\RR)=L^{u,v}\cap G(\RR)$, where $G(\RR)$ is the real
part of $G$. Consider the case  when $u=e$ and $v=w_0$, the
longest element in $W$. In this case $L^{u,v}$ is biregularly
isomorphic to the intersection of two open opposite Schubert
cells  $C_{w_0}\cap w_0C_{w_0}$, where $C_{w_0}=(Bw_0B)/B$ is the
open Schubert cell in the flag variety $G/B$. These opposite
cells appeared in the literature in various contexts (see e.g.
\cite{BFZ, R1}). Let $\sharp$ denote the number of connected
components in $L^{e,w_0}(\RR)$; following \cite{Z} we write
$\sharp=\sharp(X_n)$, where $X_n=A_n, B_n,\dots, G_2$ runs over
all types of simple Lie groups in the Cartan--Killing
classification.

The numbers $\sharp(A_n)$ were determined in \cite{SSV97, SSV98}: it
turns out that $\sharp(A_1)=2$, $\sharp(A_2)=6$, $\sharp(A_3)=20$,
$\sharp(A_4)=52$, and $\sharp(A_n)=3\cdot2^n$ for $n\gs5$.
The numbers $\sharp(D_n)$ were determined in \cite{Z}; namely,
$\sharp(D_n)=3\cdot2^n$ for $n\gs4$. It is also shown in \cite{Z} that
$\sharp(E_n)=3\cdot2^n$ for $n=6, 7, 8$. The case $G_2$ was treated in
\cite{R2}: $\sharp(G_2)=11$ (see also \cite{Z} for another proof of
this result). For nonsimply-laced series $B_n$ and $C_n$, only the
simplest case $n=2$ is known; in this case
$\sharp(B_2)(=\sharp(C_2))=8$
(see \cite{R2, Z}).

In this note we calculate $\sharp(X_n)$ for the remaining simple Lie groups
of types $B_n$, $C_n$, and $F_4$, and thus provide a complete solution for
the problem posed in \cite{Z, Remark~5.3}.

\proclaim{Theorem 1} For any  $n\gs4$ one has
$\sharp(B_n)=\sharp(C_n)=(n+5)\cdot2^{n-1}$. Besides, $\sharp(B_3)=\sharp(C_3)
=30$, $\sharp(B_2)=\sharp(C_2)=8$ and $\sharp(F_4)=80$.
\endproclaim

In fact, we prove a more general result, and find the number of
connected components of $L^{u,v}(\RR)$ for any generic pair $(u,v)\in
W\times W$ (see Theorem~4 below).

The authors would like to thank B.~Shapiro and A.~Zelevinsky for
valuable discussions and encouragement.   They express their
gratitude to Volkswagen--Stiftung for the financial support of
their stay at the Mathematischen Forschungsinstitut Oberwolfach
in Summer 2000 (under the program ``Research in Pairs'').
M.~G.~and A.~V.~are also grateful to the Gustafsson foundation
for the financial support of their visits to KTH in the Fall 2000
and in the Spring 2001.

\heading 2. Proofs \endheading

We start with reminding the  following important construction
from \cite{SSVZ, Z}; in fact, this is not the original construction
itself, but rather its version reduced modulo 2.

Let $\Pi$ be the Coxeter graph of $G$, and let $s_i (i \in
\Pi)$ be the system of simple reflections that generate $W$. A
word $\ii = (i_1, \ldots, i_m)$ in the alphabet $\Pi$ is a {\it
reduced word\/} for $w \in W$ if $w = s_{i_1} \cdots s_{i_m}$,
and $m$ is the smallest length of such a factorization. The
length of any reduced word for $w$ is called the {\it length} of
$w$ and denoted by $\l (w)$.

Let $\gg $ be the Lie algebra of $G$,
$\hh $ be the Cartan subalgebra of $\gg$  and $A = (a_{ij})$ be the
Cartan matrix. Recall that for $i \neq
j$ the indices $i$ and $j$ are adjacent in $\Pi$ if and only if
$a_{ij} a_{ji} \neq 0$; we shall denote this by $\{i,j\} \in \Pi$.

 Let us consider the group $W \times W$. It
corresponds to a graph $\tilde \Pi$ given by the union of two disconnected
copies of $\Pi$. We identify the vertex set of $\tilde \Pi$ with
$\{+1, -1\} \times \Pi$, and write a vertex $(\pm 1, i) \in \tilde
\Pi$ simply as $\pm i$. For each $i \in \Pi$, we set
$\veps (\pm i) = \pm 1$ and $|\pm i| = i$. Thus, two
vertices $i$ and $j$ of $\tilde \Pi$ are adjacent if and only if
$\veps (i) = \veps (j)$ and $\{|i|, |j|\} \in \Pi$. In this
notation, a reduced word for a pair $(u,v) \in W \times W$ is an
arbitrary shuffle of a reduced word for $u$ written in the
alphabet $-\Pi$ and a reduced word for $v$ written in the alphabet
$\Pi$. The set of all reduced words for a given pair $(u,v)\in W
\times W$ is denoted by $R(u,v)$.

Now let us fix a pair $(u,v) \in W \times W$, and let $d = \l (u)
+ \l(v)$. Let $\ii = (i_1, \ldots, i_d) \in R(u,v)$ be any
reduced word for $(u,v)$. We associate to $\ii$ an $d \times d$
matrix $(\C_{kl})$ over the two-element field $\FF_2$ in the
following way: set $\C_{kl} = 1$ if $|i_k| = |i_l|$ and $\C_{kl}
= a_{|i_k|,|i_l|}\bmod 2$ if $|i_k| \neq |i_l|$.

Next, we associate with $\ii$ a graph $\Sigma (\ii)$ on the set of
vertices $[1,d] = \{1, 2, \ldots, d\}$. For $l \in [1,d]$, we
denote by $l^- = l^-_\ii$ the maximal index $k$ such that $1 \ls
k < l$ and $|i_k| = |i_l|$; if $|i_k| \neq |i_l|$ for $1 \ls k <
l$ then we set $l^- = 0$. The edges of $\Sigma(\ii)$ are now
defined as follows.

A pair $\{k,l\} \subset [1,d]$ with $k<l$
is an edge of $\Sigma(\ii)$ if it satisfies one of the following
three conditions:

(i) $k=l^-$;

(ii) $k^- < l^- < k$, $\{|i_k|, |i_l|\} \in \Pi$, and
$\veps(i_{l^-})=\veps(i_{k})$;

(iii) $l^-< k^- < k$, $\{|i_k|, |i_l|\} \in \Pi$, and
$\veps(i_{k^-})=-\veps(i_{k})$.

\noindent The edges of type (i) are called {\it horizontal}, and
those of types (ii) and (iii) {\it inclined}. Each inclined edge
corresponds to an edge of the graph $\Pi$.
We shall write $\{k, l\} \in \Sigma(\ii)$ if $\{k, l\}$ is an
edge of $\Sigma(\ii)$.

We now associate to each $r \in [1,d]$ a transvection
$\tau_r \: \FF_2^d \to \FF_2^d$
defined as follows: $\tau_r
(\xi_1, \dots, \xi_d) = (\xi'_1, \dots, \xi'_d)$, where $\xi'_k =
\xi_k$ for $k \neq r$, and
$$
\xi'_r  = \xi_r  + \sum_{\{k, r\} \in
\Sigma(\ii)}\C_{kr} \xi_k. \tag1
$$
(note that (1) coincides with the reduction modulo 2 of
formula (2.2) in \cite{Z}).
We call an index $r \in
[1,d]$ {\it $\ii$-bounded\/} if $r^- > 0$.
The set of all bounded indices (and corresponding vertices of
$\Sigma(\ii)$) is denoted by $B$ and its complement is denoted by $C$.

Let $\Gamma_{\ii}$
denote the group of linear transformations of $\FF_2^d$ generated by
the transvections $\tau_r$ for all $\ii$-bounded indices $r \in
[1,d]$. The following result was conjectured in \cite{SSVZ} for a simply
laced case, and proved in \cite{Z} in the general case (see also
\cite{SSV97} for the case of open cells for type $A_n$).

\proclaim{Theorem 2}
For every reduced word $\ii \in R(u,v)$,
the connected components of $L^{u,v}(\RR)$ are in a natural bijection
with the $\Gamma_\ii $-orbits in $\FF_2^{d}$.
\endproclaim

This theorem,   together
with the description of orbits of groups generated by symplectic
transvections presented in \cite{SSV98, SSVZ},  form the basis of the
enumerative results in the simply-laced case cited in the Introduction.

However, in the nonsimply-laced case, the transvections generating
$\Gamma_{\ii}$ are no longer symplectic. To handle this case, we
have to extend several results of \cite{SSV98, SSVZ}.

Let $W^t$, $t\in [1,n-1]$, be the Coxeter group with $n$
generators $s_1,\dots,s_n$ and relations of the form $s^2_i=1$,
$(s_is_j)^2=1$ for $j>i+1$, $(s_is_{i+1})^3=1$ for $i\ne t$, and
$(s_ts_{t+1})^4=1$. Denote by $\Pi^t$ the Coxeter graph of $W^t$.
Finally, define the $n\times n$ matrix $A^t$ as follows:
$a_{t,t+1}=2$, $a_{t+1,t}=1$, $a_{i,i+1}= a_{i+1,i}=1$ for any
$i\in [1, n-1]$, $i\ne t$, $a_{ij}=0$ for $|i-j|\ne1$.

Fix a pair $(u,v)\in W^t\times W^t$,
take an arbitrary reduced word $\ii$ for the pair $(u,v)$, and build the
graph $\Sigma(\ii)$ and transvections $\tau_r$
exactly as above, with $\Pi$ replaced by $\Pi^t$ and $A$ replaced by $A^t$.
Observe that for $t=1$ the above construction describes
the $C_n$ case, for $t=n-1$, it describes the $B_n$ case, and for $n=4$, $t=2$,
the $F_4$ case.

Define $\Pi_U^t$  to be the subgraph of $\Pi^t$ induced by the
vertices $\{1,2,\dots,t\}$, and $\Pi_L^t$ to be the complement to
$\Pi_U^t$ in $\Pi^t$. In accordance with this partition of
$\Pi^t$, we subdivide the vertex set of $\Sigma$ into $U=\{k\in
\Sigma \: |i_k|\in\Pi_U^t\}$ and its complement $L$ (we omit in
the notation the dependence of $\Sigma$ and other objects on the
reduced word $\ii$ which is assumed fixed). Together with the
partition into bounded and unbounded vertices described above,
this gives four subsets, which we denote $B_U$, $B_L$, $C_U$, and
$C_L$; the subgraph of $\Sigma$ induced by a subset $X\subseteq
\Sigma$ is denoted $\Sigma_X$, and $\FF_2^X$ is the linear
subspace of $\FF_2^d$ defined by the condition that all
coordinates that correspond to $\Sigma \setminus X$ vanish.
 The subgroups $\Gamma_U$ and $\Gamma_L$
of $\Gamma$ are defined in a natural way; clearly, $\Gamma$ is
generated by $\Gamma_U$ and $\Gamma_L$.

For any vector $\nu\in\FF_2^L$, the action of $\Gamma_U$ preserves
the affine subspace $\nu+\FF_2^U$. Identifying $\nu+\FF_2^U$ with
$\FF_2^U$ with the help of the shift $\xi\mapsto\xi-\nu$, we get
an action of $\Gamma_U$ on $\FF_2^U$; slightly abusing notation,
we call it the {\it $\Gamma_U(\nu)$-action on\/} $\FF_2^U$. Note
that for $\nu\ne0$ the $\Gamma_U(\nu)$-action is not linear, but
rather affine; the $\Gamma_U(0)$-action coincides with the usual
linear action of $\Gamma_U$ on $\FF_2^U$.

It follows from \cite{SSVZ, Proposition~6.1} that the number of
fixed points of the $\Gamma_U(0)$-action equals $2^t$; the number
of nontrivial orbits of this action (those which are not fixed
points) we denote by $N_U$. In a similar fashion, define the
number $N_L$ of nontrivial orbits of the action of $\Gamma_L$ on
$\FF_2^L$; the number of fixed points of this action equals
$2^{n-t}$. Observe that one can define also the
$\Gamma_L(\varkappa)$-action on $\FF_2^L$ for any
$\varkappa\in\FF_2^U$, but this action does not depend on
$\varkappa$ and coincides with the $\Gamma_L$-action.

\proclaim{Lemma 1} For any vector $\nu\in\FF_2^L$ there are
$2^t+N_U$ orbits of the $\Gamma_U(\nu)$-action on $\FF_2^U$.
\endproclaim

\demo{Proof} Indeed, the $\Gamma_U(\nu)$-action on $\FF_2^U$ is
generated by affine transformations of the form
$\theta_j(\xi)=\tau^U_j(\xi)+b_j$ for $j\in B_U$, $\xi\in\FF_2^U$,
where $\tau^U_j$ is the symplectic transvection with respect to
the restriction of $\Omega$ to $\FF_2^U$ and $b_j$ depends only
on $\nu$. Assume that $\xi^*\in\FF_2^U$ is a fixed point of this
affine action. Then $\tau^U_j(\xi)-\xi^*= \tau^U_j(\xi-\xi^*)$
for any $j\in B_U$ and $\xi\in \FF_2^U$, and hence the orbits of
the $\Gamma_U(\nu)$-action  are just the orbits of the
$\Gamma_U(0)$-action shifted by $\xi^*$. Therefore, the number of
affine orbits equals $2^t+N_U$, the number of
$\Gamma_U(0)$-orbits.

It remains to check the existence of a fixed point of the affine
action. Such a fixed point should satisfy the equation
$M\xi=b(\nu)$ for some $b(\nu)\in\FF_2^{B_U}$, where $M\:\FF_2^U\to
\FF_2^{B_U}$ is given by
$$
(M\xi)_j= \xi_j  + \sum_{\{k, j\} \in\Sigma_U}\xi_k.
$$
The kernel of $M$ consists of the fixed points of the
$\Gamma_U(0)$-action. Therefore, its dimension equals to
$t=|C_U|$, which means that the image of $M$ coincides with
$\FF_2^{B_U}$. Therefore, equation $M\xi=b$ can be solved for any
$b$, and we are done.
\qed
\enddemo

The number of $\Gamma$-orbits in $\FF_2^d$ is determined as follows.

 \proclaim{Theorem 3} Assume that $\Sigma_B$ is connected,
then the number of $\Gamma$-orbits in $\FF_2^d$ equals
$2^n+2^{n-t}N_U+2^tN_L$.
\endproclaim

\demo{Proof}  Observe first that the intersections of the orbits of
$\Gamma$-action with $\FF_2^L$ are exactly the orbits of the
$\Gamma_L$-action. Consider first $\Gamma$-orbits whose intersections with
$\FF_2^L$ are fixed points of the $\Gamma_L$-action. The number
of fixed points of the $\Gamma_L$-action is $2^{n-t}$, hence by Lemma 1
we see that the number of such $\Gamma$-orbits equals to
$(2^t+N_U)2^{n-t}$.

Next, consider $\Gamma$-orbits whose intersections with $\FF_2^L$
are nontrivial $\Gamma_L$-orbits. We claim that the number of
such $\Gamma$-orbits equals $2^{t}N_L$.

Indeed, let us fix a vector $\nu\in\FF_2^L$ in such a
$\Gamma_L$-orbit, and consider the $\Gamma_U(\nu)$-action on
$\FF_2^U$. As before, by Lemma 1 we get an affine action having
$2^t+N_U$ orbits for this choice of $\nu$. We shall show that the
$\Gamma_L$-action can be used to glue these orbits into $2^t$
$\Gamma$-orbits differing only by the values on $C_U$.
To achieve this, it is enough to show that one can change the value
$\xi_r$ for any given $r\in B_U$, and to keep all the other
$\xi_j$, $j\in B$, unchanged. This is evidently true if $\tau_r(\xi)\ne\xi$,
so in what follows we assume that $\tau_r(\xi)=\xi$.

Denote by $T(\xi)$ the set of all $j\in B_L$ such that
$\tau_j(\xi)\ne\xi$; $T(\xi)\ne\varnothing$, since $\nu$ belongs
to a nontrivial $\Gamma_L$-orbit. The connectivity of $\Sigma_B$
implies the existence of a path joining $r$ with the set
$T(\xi)$. Moreover, since the set of all vertices $q$ having the
same height $|i_q|$ is connected in $\Sigma_B$, there exists a
{\it monotone\/} path from $r$ to $T(\xi)$, that is, a one for
which the height changes monotonously along the path. Let
$P=(q_0\in T(\xi),q_1,\dots,q_k=r)$ be a shortest monotone path
between $T(\xi)$ and $r$; besides, let $q_l$ be the first vertex
at height $t$ in this path. Note that since the path $P$ is
monotone, all the vertices $q_j$, $j\in [l,k]$, belong to $B_U$.

Assume first that $\tau_{q_j}(\xi)=\xi$ for $j\in [l,k]$. Apply
consequently $\tau_{q_0},\tau_{q_1},\dots,\tau_{q_k}$; upon
applying $\tau_{q_i}$, the value $\xi_{q_i}$ is changed, since
the only edge of the type $\{q_i,q_j\}$, $j<i$, is the edge
$\{q_i,q_{i-1}\}$ (otherwise the path is not the shortest
possible). Hence, applying the whole sequence results in changing
the value $\xi_r$. To restore the values $\xi_{q_i}$, $i\ne
[0,k-1]$, apply consequently
$\tau_{q_{l-1}},\tau_{q_{l-2}},\dots, \tau_{q_0}$ followed by
$\tau_{q_l},\tau_{q_{l+1}},\dots,\tau_{q_{k-1}}$.

Otherwise, let $q_m$, $m\in [l,k]$, be the vertex of $P$ closest to $r$ for
which $\tau_{q_m}(\xi)\ne\xi$. Apply consequently $\tau_{q_m},\tau_{q_{m+1}},
\dots,\tau_{q_k}$ to change the value $\xi_r$. To restore the values
$\xi_{q_i}$, $i\in [m,k-1]$, we have
to solve the same problem as above, but now the length of a shortest monotone
path to $T(\xi)$ equals $k-1$, and we are done by induction.

Proceeding in this way, we see that any $\Gamma$-orbit whose
intersection with $\FF_2^L$ does not coincide with a fixed point of the
$\Gamma_L$-action contains a vector that vanishes at any point of
$B_U$. Therefore, the only invariants of such an orbit are the values of
$\xi$
at the points of $C_U$. Since the number of these points equals $t$,
we get $2^t$ $\Gamma$-orbits per each nontrivial $\Gamma_L$-orbit,
which totals to $2^tN_L$ $\Gamma$-orbits.
\qed
\enddemo

To prove our main result, we need a description of the action of $\Gamma$
in the case when $\Sigma$ is a path. So, let $\Sigma=P$ be a path on $m$
vertices $\{1,2,\dots,m\}$. The action of $\Gamma$ on $\FF_2^P=\FF_2^m$ is
generated by symplectic  transvections
 $\tau^P_j$, $j\in [2,m]$, defined by
$$
\tau^P_j (\zeta) = \zeta +
(\zeta_{j-1} + \zeta_{j+1}) e_j,
$$
where $\{e_j\}$ is the standard basis of $\FF_2^m$; we call it the {\it
$\Gamma_P$-action}.

\proclaim{Lemma 2} The number of orbits of the $\Gamma_P$-action equals
$m+1$. Exactly two of these orbits are fixed points of the $\Gamma_P$-action.
\endproclaim

\demo{Proof} Let $\zeta\in \FF_2^m$ be of the form
$$
\zeta=(
\underbrace{0\ldots 0}_{l_1} \underbrace{1\ldots 1}_{m_1}\ldots
\underbrace{0\ldots 0}_{l_k} \underbrace{1\ldots
1}_{m_k}\underbrace{0\ldots 0}_{l_{k+1}}),
$$
where $l_1,l_{k+1}\gs 0$, $l_2,\ldots, l_k, m_1,\ldots,m_k > 0$;
we put $c(\zeta)=k$. It is easy to see that $c(\tau^P_j (\zeta))=c(\zeta)$,
and that $(\tau^P_j(\zeta))_1=\zeta_1$. Let us prove that if
$c(\zeta)=k$ and $\zeta_1=1$ (resp.,~$\zeta_1=0$),
then there exists $\gamma\in \Gamma_P$ such that
$\gamma(\zeta)=(\underbrace{1,0,1,0,\ldots,1,0,1}_{2k-1},0,\allowmathbreak
\ldots,0)$
(resp.,~$\gamma(\zeta)=(\underbrace{0,1,0,1,\ldots,0,1}_{2k},0,\ldots,0)$).

Indeed, if $\zeta=(\zeta_1,\ldots,\zeta_{j-1},\underbrace{1,\ldots,1}_{l},0,
\zeta_{j+l+1},\ldots)$, then
$$
\tau^P_{j+1}\cdots\tau^P_{j+l-1}(\zeta)=
(\zeta_1,\ldots,\zeta_{j-1},1,\underbrace{0,\ldots,0}_{l},
\zeta_{j+l+1},\ldots).
$$
Similarly, if
$\zeta=(\zeta_1,\ldots,\zeta_{j-1},\underbrace{0,\ldots,0}_{l},1,
\zeta_{j+l+1},\ldots)$,
then
$$
(\tau^P_{j+2}\tau^P_{j+1})\cdots(\tau^P_{j+l}\tau^P_{j+l-1})(\zeta)=
(\zeta_1,\ldots,\zeta_{j-1},0,1,\underbrace{0,\ldots,0}_{l-1},
\zeta_{j+l+1},\ldots).
$$
Combining transformations of these two types, we can eventually
bring $\zeta$ to the required form.

Since the number of these forms equals $m+1$, and any two of them differ
either at $c(\zeta)$, or at $\zeta_1$ (or at both of them), we conclude that
the number of $\Gamma_P$-orbits equals $m+1$. Evidently, if $m$ is even,
then $(0,1,\dots,0,1)$ is a fixed point of the $\Gamma_P$-action, while
if $m$ is odd, the $(1,0,\dots,1,0,1)$ is such a fixed point. The only
other fixed point is $(0,\dots,0)$.
\qed
\enddemo

Now we return to the cases $B_n$ and $C_n$. We say that a pair
$(u,v)$ is {\it generic\/} if there exists $\ii\in R(u,v)$ such
that the subgraph  $\Sigma_B(\ii)$ is connected, and the subgraph
$\Sigma_{B_L}(\ii)$ (in the $C_n$ case) or $\Sigma_{B_U}(\ii)$
(in the $B_n$ case) is $E_6$-compatible. One can prove easily
that almost all pairs $(u,v)$ are generic (cp.~with the similar
result in the $A_n$-case proved in \cite{SSV99}). Recall that in
the $C_n$ (respectively, $B_n$) case, the graph $\Sigma_U$
(respectively, $\Sigma_L$) is a path. Let $m$ denote the number
of vertices in $U$ for the $C_n$ case, and the number of vertices
in $L$ for the $B_n$ case. It is easy to see that this number
depends only on the pair $(u,v)$, and does not depend on the
reduced word $\ii\in R(u,v)$.

\proclaim{Theorem 4} Let $(u,v)$ be a generic pair, then the number of
connected components in $L^{u,v}(\RR)$ equals $(m+5)\cdot2^{n-1}$ for both
types $B_n$ and $C_n$.
\endproclaim

\demo{Proof} Since the pair $(u,v)$ is generic, there exists $\ii
\in R(u,v)$ such that the subgraph of $\Sigma(\ii)$ induced by
$B$ is connected. Hence, by Theorem~3, the number of
$\Gamma_\ii$-orbits for type $C_n$ equals $2^n+2^{n-1}N_U+2N_L$,
and for type $B_n$, equals $2^n+2^{n-1}N_L+2N_U$. Besides, by
\cite{SSVZ, Th.~7.2}, the $E_6$-compatibility condition implies
that $N_L=2^n$ for type $C_n$, and $N_U=2^n$ for type $B_n$.
Moreover, by Lemma~2, $N_U=m-1$ for type $C_n$, and $N_L=m-1$ for
type $B_n$. Therefore, in both cases the total number of orbits
equals $2^n+(m-1)\cdot2^{n-1}+2^{n+1}= (m+5)\cdot2^{n-1}$. By
Theorem~2, this number equals the number of connected components
in $L^{u,v}(\RR)$. \qed
\enddemo

To prove Theorem~1 stated in the introduction one has to check
that the pair $(e,w_0)$ is generic for $n\gs4$. This fact follows
immediately from Figure~1 presenting the graph $\Sigma(\ii)$ and
its corresponding subgraphs for $n=4$ and $\ii=1234123412341234$.

Consider now the cases $n=2,3$. One can check easily that the subgraphs
$\Sigma_B$ remain connected, though the pair $(e,w_0)$
is no longer generic; therefore, Theorem~3 remains valid.
Besides, one gets $N_U=N_L=1$ for types $B_2$ and $C_2$,
$N_U=2$, $N_L=7$ for type $C_3$, and $N_U=7$, $N_L=2$ for type $B_3$.
Thus, Theorem~3 yields $\sharp_2=4+2+2=8$ and $\sharp_3=8+8+14=30$.

\vskip 10pt
\centerline{\hbox{\epsfxsize=11cm\epsfbox{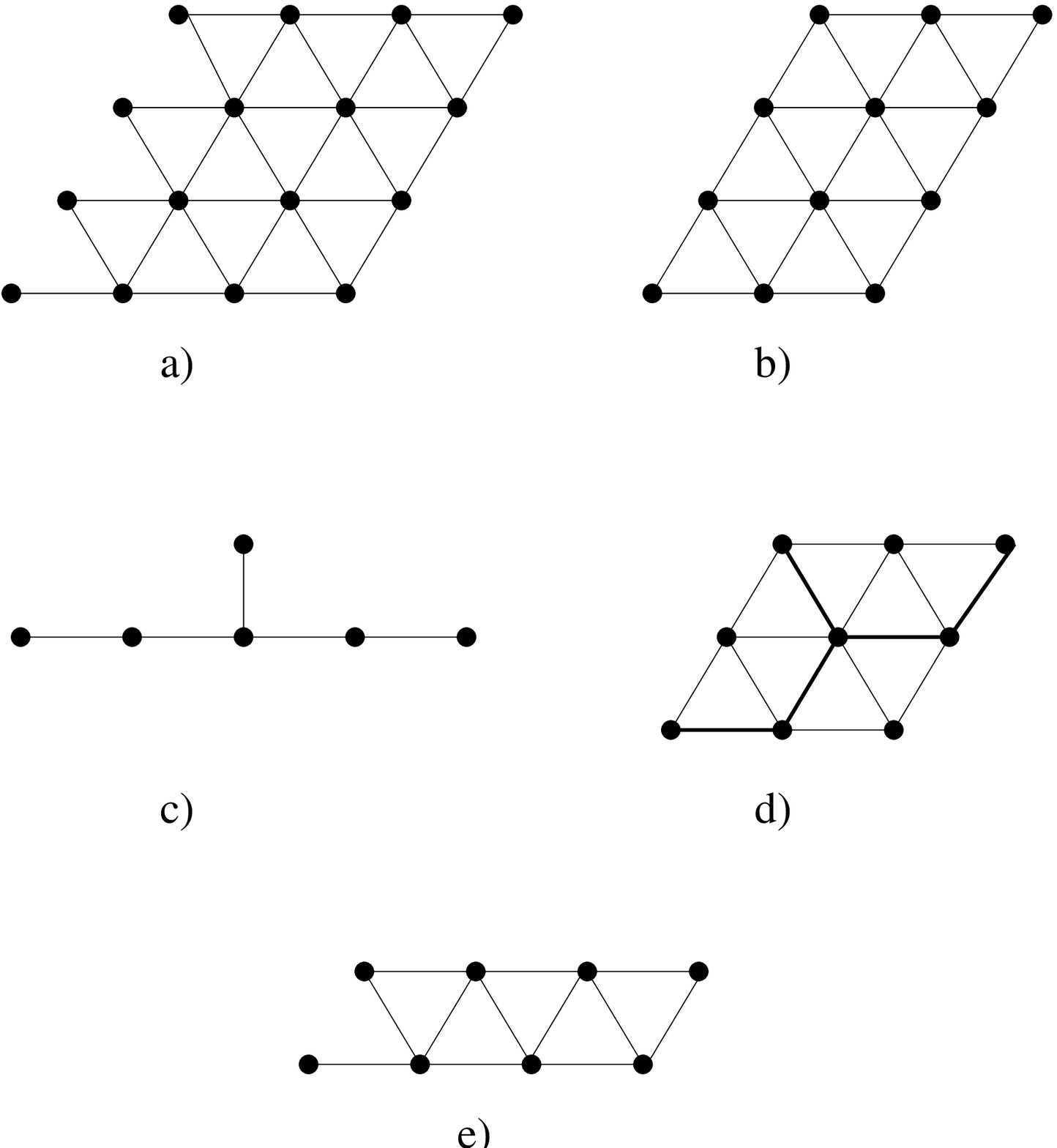}}}
\midspace{1mm}
\caption{Fig.~1. To the proof of Theorem~1: {\rm a)} graph $\Sigma(\ii)$;
{\rm b)}~graph $\Sigma_B(\ii)$; {\rm c)} graph $E_6$;
{\rm d)}~graph $\Sigma_{B_L}(\ii)$ and an induced $E_6$ in it;
{\rm e)}~graph $S(4)$}
\vskip 5pt

To treat the case of $F_4$, we first  consider the a
graph $S=S(m)$ defined as follows: $S$ contains vertices
$\{1,\ldots,2m\}$ arranged into two levels, the lower (resp.
upper) level is formed by odd-numbered (resp. even-numbered)
vertices. Horizontal edges are of the form $(2i, 2i+2)$ and
$(2i-1, 2i+1)$, and inclined edges are of the form $(2i+1, 2i+2)$
and $(2i, 2i+1)$, where $i$ runs from $1$ to $m-1$ (see Figure~1e).
It is convenient to represent elements of $ \FF_2^{S}$
by vectors $\zeta=(\zeta_{i})_{i=1}^{2m}\in \FF_2^{2m}$. The
$\Gamma_S$-action on $ \FF_2^{S}$ is generated
by transvections $\tau^{S}_j$, $j\in [3,2m]$, defined by
$$
\tau^{S}_{j} (\zeta) = \zeta +
(\zeta_{j-2} + \zeta_{j-1}+\zeta_{j+1} + \zeta_{j+2} ) e_{j}\ . \tag2
$$
where $\zeta_i=0$ if $i>2m$.

\proclaim{Lemma 3} Let $m > 2$, then every nontrivial orbit of
the $\Gamma_S$-action
contains either an element of the form
$(\zeta_1,\zeta_2,\zeta_3,\zeta_4,0,\dots,0)$
where not all $\zeta_i$ are equal to zero,
or the element $\bar\zeta=(0,0,1,1,1,0,\ldots,0)$.
\endproclaim

\demo{Proof} Let us fix a nontrivial orbit $\OO$ of the $\Gamma_S$-action.
To prove the statement, it suffices to show that for any $\xi\in \OO$ of the
form
$$\xi= (\xi_1,\ldots,\xi_{j-1},1,\underbrace{0\ldots 0}_{2m-j})$$
such that $j>4$ and $\xi
\ne \bar\zeta$, there exists $\gamma\in\Gamma_S$ such that $\gamma(\xi)_i=0$
for $i\gs j$.

If the set $T=\{ i\: 3\ls i\ls j, \tau^{S}_{i}(\xi)\ne \xi \}$ is not empty
(this is clearly the case for $j=2m$), we denote by $k$ the largest element
in $T$ and define $\gamma$ as the product of $\tau_i^S$ along any shortest
path from $k$ to $j$.
Then $\gamma(\xi)_i=0$ for $i\gs j$.

Otherwise, $T=\varnothing$ and the smallest $i$ such that
$\tau^{S}_{i}(\xi)\ne \xi $
is equal either to $j+1$ or to $j+2$. In the first case, $\xi$ has to be of
the form
$$\xi= (\xi_1,\ldots,\xi_{j-5},0,1,0,0,1,0,\underbrace{0\ldots 0}_{2m-j-1}).$$
Define
$\gamma=\tau^{S}_{j+1}\tau^{S}_{j}\tau^{S}_{j-2}\tau^{S}_{j-1}\tau^{S}_{j+1}$. Then
$$\gamma(\xi)= (\xi_1,\ldots,\xi_{j-5},
0,1,1,1,\underbrace{0\ldots 0}_{2m-j+1})\ .$$
In the second case, either $\xi=\bar\zeta$ and we are done, or
$$\xi= (\xi_1,\ldots,\xi_{j-6},0,0,0,1,1,1,0,0,
\underbrace{0\ldots 0}_{2m-j-2}),$$
in which case we
put $\gamma=\tau^{S}_{j+1}\tau^{S}_{j-1}\tau^{S}_{j-3}\tau^{S}_{j+2}\tau^{S}_{j}\tau^{S}_{j-1}
\tau^{S}_{j+1}\tau^{S}_{j+2}$.
Then
$$\gamma(\xi)= (\xi_1,\ldots,\xi_{j-5},0,0,1,1,1,
\underbrace{0\ldots 0}_{2m-j+1}).$$
This finishes the proof.
\qed
\enddemo

\proclaim{Corollary} If $m>2$ then the number of orbits of the
$\Gamma_S$-action is equal to $12$. Four of these orbits are fixed points of
the action.
\endproclaim

\demo{Proof} It follows from (2) that for every choice of
$\alpha,\beta \in \FF_2$ there is
exactly one fixed point of the $\Gamma_S$-action with $\zeta_{2n-1}=\alpha$,
$\zeta_{2n}=\beta$.
Thus, we have four orbits that are fixed points of the action.

By the previous lemma, any other orbit is either the orbit through
$\bar\zeta$ or the orbit through
an element of the form $\zeta=(\zeta_1,\zeta_2,\zeta_3,\zeta_4,0,\ldots,0)$,
where $\zeta_i$ cannot be all equal to zero. It is easy to see that if
$\zeta_3\ne\zeta_2$ then  $\tau_4^S(\zeta)\ne\zeta$; moreover, either
$\tau_3^S(\zeta)\ne\zeta$, or $\tau_3^S\tau_4^S(\zeta)\ne\tau_4^S(\zeta)$.
Besides, if $\zeta_3=\zeta_2$ and $\zeta_4=\zeta_1+\zeta_2$ then
$\tau_3^S(\zeta)=\tau_4^S(\zeta)=\zeta$.
This means that the number of nontrivial orbits
does not exceed $8$.

A non-homogeneous quadratic form
$$Q_S(\xi)=\sum_{i\in S} \xi_i+\sum_{(i,j)\in S} \xi_i \xi_j$$
is an invariant of the $\Gamma_{S}$-action (see \cite{SSVZ}),
along with the values of $\xi_1, \xi_2$.
Now, to finish the proof it is sufficient to notice that the triple
$(\xi_1, \xi_2, Q_S(\xi))$
takes different values on the following eight elements:
$$ \align &(1,1,1,1,0,\dots,0),  (1,1,1,0,0,\dots,0),
           (1,0,0,1,0,\dots,0),  (1,0,0,0,0,\dots,0), \\
          &(0,1,1,1,0,\dots,0),  (0,1,1,0,0,\dots,0),
           (0,0,1,0,0,\dots,0),  (0,0,1,1,1,0,\dots,0).
\endalign
$$
\qed
\enddemo

We are now in a position to finish the proof of ~Theorem 1.

\proclaim{Theorem 5} $\sharp(F_4)=80$.
\endproclaim

\demo{Proof} Recall, that $\ii=(1234)^6$ is a reduced word for $w_0$ in the
Weyl
group that correspond to $F_4$. We can use Theorem 3 again. In this case,
$n=4$, $t=2$
and both subgraphs $\Sigma_L$ and $\Sigma_U$ coincide with $S(6)$. Then,
by Theorem 3 and Corollary to Lemma 3, $\sharp(F_4)= 2^4 + 2*2^2*8=80$.
\qed
\enddemo

\enddocument